\pdfoutput=1

\documentclass[12pt]{amsart}
\usepackage[centertags]{amsmath}
\usepackage{amsfonts}
\usepackage{amsthm}

\usepackage{amscd}
\usepackage{amssymb}

\theoremstyle{plain}
\newtheorem{thm}{Theorem}[section]
\newtheorem{cor}[thm]{Corollary}
\newtheorem{lemma}[thm]{Lemma}
\newtheorem{prop}[thm]{Proposition}
\theoremstyle{definition}

\newtheorem{rem}[thm]{Remark}
\newtheorem{ex}[thm]{Example}
\let\oldrem\rem
\renewcommand{\rem}{\oldrem\normalfont}

\theoremstyle{definition}

\newtheorem{remark}[thm]{Remark}

\numberwithin{equation}{section}

\DeclareMathOperator{\Hom}{Hom}
\DeclareMathOperator{\HH}{H}

\newcommand\TT{{\mathbb{T}}}
\newcommand\PP{{\mathbb{P}}}

\newcommand\QQ{{\mathbb{Q}}}
\newcommand\ZZ{{\mathbb{Z}}}
\newcommand\CC{{\mathbb{C}}}
\newcommand\BB{{\mathbb{B}}}
\newcommand\EE{{\mathbb{E}}}

\newcommand\cX{{\mathcal X}}
\newcommand\cY{{\mathcal Y}}
\newcommand\cO{{\mathcal O}}

\newcommand\cN{{\mathcal N}}
\newcommand\cC{{\mathcal C}}
\newcommand\cT{{\mathcal T}}

\newcommand\ra{{\ \rightarrow\ }}
\newcommand\lra{\longrightarrow}
\newcommand\lla{\longleftarrow}
\newcommand\iso{{\ \cong\ }}

\def\min{\text{\rm min}}

\begin{document}

\title{On rigidity of flag varieties}
\author[A.~Weber]{Andrzej Weber}
\author[J.~A.~Wi\'{s}niewski]{Jaros\l{}aw A.~Wi\'{s}niewski}
\address{University of Warsaw, Institute of Mathematics,
Banacha 2, Warszawa, Poland \\ Institute of Mathematics, Polish Academy of Sciences,  \'Sniadeckich 8, Warszawa, Poland}
\email{A.Weber@mimuw.edu.pl}
\address{University of Warsaw, Institute of Mathematics,
Banacha 2, Warszawa, Poland}
\email{J.Wisniewski@mimuw.edu.pl}
\keywords{}
\subjclass[2000]{}

\begin{abstract}
  We prove that the variety of complete flags for any semisimple
  algebraic group is rigid in any smooth family of Fano
  manifolds.
\end{abstract}

\maketitle

\section{Results}

The aim of this note is to prove the following rigidity theorem.
\begin{thm}\label{rigidity}
  Let $G$ be a semisimple algebraic group defined over complex
  numbers. By $B<G$ we denote its Borel subgroup. Let
  $\pi: \cX\ra\Delta\ni 0$ be a smooth family of complex projective
  manifolds over a positive-dimensional connected base $\Delta$. That
  is $\pi:\cX\to\Delta$ is a submersion of smooth varieties and for
  every $t\in\Delta$ we set $X_t=\pi^{-1}(t)$. Assume that the
  relative anticanonical divisor $-K_{\cX/\Delta}$ is $\pi$-ample and
  for $t\ne 0$ the variety $X_t=\pi^{-1}(t)$ is isomorphic to the
  variety of complete flags $G/B$. Then $X_0\iso G/B$
\end{thm}

The rigidity of rational homogeneous varieties of type
$G/P_\text{max}$, where $P_\text{max}<G$ is a maximal parabolic
subgroup, was studied by Siu, Hwang and Mok, see \cite{SiuCrelle},
\cite{HwangMokInvMath}, \cite{HwangMokAnnSciENS}. Their results
concern irreducible Hermitian symmetric spaces. As it was observed by
Pasquier and Perrin in \cite[Prop.~2.3]{PasquierPerrin} not all
varieties of this type are rigid. Also, if a parabolic subgroup $P<G$
is neither Borel nor maximal then $G/P$ may not be rigid. Moreover,
the theorem fails to be true if we do not assume ampleness of the
relative anticanonical divisor. We discuss these issues in the last
section of the present paper.

The proof of the theorem is based on the characterization of varieties
of complete flags by Occhetta, Sola Conde, Watanabe and the second
author, from \cite{OSWW}, and the rigidity of the nef cone for Fano
manifolds from \cite{BLMS}. The main theorem of \cite{OSWW}
characterizes flag manifolds of type $G/B$ as the unique Fano
manifolds whose all elementary contractions in terms of Mori theory
are $\PP^1$-bundles, see theorem \ref{charOSWW}.  The main technical
result of the paper is the following equidimensional rigidity for
Fano-Mori fiber bundles.

\begin{prop}\label{equidimensional}
  Let $\pi: \cX\ra\Delta\ni 0$ be a smooth family of complex
  projective manifolds over a positive-dimensional connected base
  $\Delta$.  We set $X_t=\pi^{-1}(t)$. Let $\cY$ be a normal variety
  with a projective morphism $\cY\ra\Delta$, we denote
  $Y_t=\pi^{-1}(t)$. Let $\varphi: \cX\ra\cY$ be a projective morphism
  of varieties over $\Delta$, which commutes with the morphisms
  $\cX\ra\Delta$ and $\cY\ra\Delta$. By $\varphi_t: X_t\ra Y_t$ we
  denote the restriction of $\varphi$. We assume the following:
  \begin{enumerate}
  \item the relative anti-canonical divisor $-K_{\cX/\Delta}$ is
    $\varphi$-ample,
  \item for every $t\ne 0$ the morphism $\varphi_t: X_t\ra Y_t$ is a
    fiber bundle, such that the fundamental group of $Y_t$ acts trivially on the
    cohomology of the fiber,
  \item for every $t\ne 0$ the rational cohomology ring $\HH^*(Y_t)$
  is generated by its first and second gradation.
  \end{enumerate}
  Then the morphism $\varphi$ is equidimensional, that is every fiber
  of $\varphi$ is of the same dimension.
\end{prop}
As a consequence of Fujita's theorem, \cite[Thm.~2.12]{Fujita} we get
the following.
\begin{cor}\label{proj-bundle}
  In the situation of proposition \ref{equidimensional} assume that
  for $t\ne 0$ the morphism $\varphi_t$ is a Zariski projective
  bundle. Then $\varphi_0$ is of the same type as well.
\end{cor}
Here a Zariski projective bundle is understood as projectivization of
a vector bundle.

We note that the condition that rational cohomology ring $\HH^*(Y_t)$
is generated by its first and second gradation is very
strong. However, it is known to be satisfied by abelian varieties,
toric varieties and complete flags, among others.  In order to prove the
rigidity theorem for complete flags we need the following.

\begin{thm}\label{generation}
  Let $G$ be a semisimple algebraic group over complex numbers and
  $P_{\min}<G$ a minimal parabolic subgroup of $G$ which contains a
  Borel group of $G$ properly. Then the rational cohomology ring
  $\HH^*(G/P_{\min})$ is generated as $\QQ$-algebra by its second
  gradation, that is $\HH^2(G/P_{\min})$.
\end{thm}

\subsection*{Notation} We deal with varieties defined over complex
numbers.  By $\iso$ we denote their isomorphism while $\simeq$ stands
for homeomorphism.  Unless specified otherwise, the cohomology is
considered with rational coefficients; given $X$, a projective
(compact) variety, by $\HH^*(X)$ we denote its (graded) cohomology
ring with $\smile$-product. A fiber bundle means a locally trivial
bundle in the underlying classical topology.

\subsection*{Acknowledgments}
The authors were supported by the research grant from Polish National
Science Center, number 2013/08/A/ST1/00804. The research was done
during miniPAGES activities at Banach Center supported by grant 346300
for IMPAN from the Simons Foundation and the matching 2015-2019 Polish
MNiSW funds. Our interest in the topic was triggered by a talk by
J.-M. Hwang during these activities who also informed us about
reference \cite{PasquierPerrin}.

\section{Proofs}
In the proofs we use the notation introduced in the previous section.

\subsection{Proof of theorem \ref{rigidity}}
Let us recall basic facts of Mori theory and refer the reader to more
comprehensive sources, e.g.~\cite{KollarMori}, for details. A
contraction $\varphi: X\ra Y$ is a surjective morphism of normal
projective varieties with connected fibers. We will also consider a
relative case when both varieties are not necessarily projective but
they admit projective morphisms $X\ra\Delta$, $Y\ra\Delta$ and
$\varphi$ commutes with morphisms to $\Delta$. The variety $X$ will be
assumed smooth (over $\Delta$) and its anti-canonical divisor $-K_X$
or, in the relative case, $-K_{X/\Delta}$ will be assumed to be
$\varphi$-ample. Then $\varphi$ is called Fano-Mori contraction. We
will say that $\varphi$ is elementary, or extremal, if
$Pic(X/Y)\iso\ZZ$. This is equivalent to say that cohomology classes
of curves contracted by $\varphi$ to points span 1-dimensional
subspace in $\HH^{2\dim X-2}(X)$.

A smooth projective variety $X$ is Fano if $-K_X$ is ample. Mori cone
theorem asserts that if $X$ is Fano then the convex cone $\cC(X)$
spanned in $\HH^{2n-2}(X)$ by the classes of curves is rational
polyhedral.  Kawamata-Shokurov contraction theorem asserts that in
this case there is a bijection between faces of $\cC(X)$ and
contractions of $X$: given a face $\Phi\subseteq\cC(X)$ there exists a
contraction $\varphi_\Phi: X\ra Y_\Phi$ which contracts to points
exactly these curves whose classes are in $\Phi$. Elementary
contractions of $X$ are associated to 1-dimensional faces of $\cC(X)$
(extremal rays).

We note that everything which is said in the previous paragraph holds
in the relative case too. In particular, if $\varphi: X\rightarrow Y$
is a contraction of a smooth variety and $-K_X$ is $\varphi$-ample, then
the cone $\cC(X/Y)$ spanned by the classes of curves contracted by
$\varphi$ is rational polyhedral and its faces are in bijection with
contractions factoring $\varphi$.

The main ingredient of the proof of \ref{rigidity} is the following
main result of \cite{OSWW}.
\begin{thm}\label{charOSWW}
  Let $X$ be a Fano manifold such that every elementary contraction of
  $X$ is a smooth $\PP^1$-fibration. Then $X$ is isomorphic to a
  complete flag variety $G/B$, where $G$ is a semisimple algebraic
  group and $B$ a Borel subgroup.
\end{thm}
It is well known that if $X\iso G/B$ then $X$ is Fano and every
elementary contraction of $X$ is a $\PP^1$-bundle $G/B\ra G/P_{\min}$,
where $P_{\min}$ is a minimal parabolic subgroup of $G$ which
contains $B$ properly. Given a semisimple group $G$, the homogeneous
spaces $G/P$ are distinguished by a set of nodes in the Dynkin diagram
of $G$. In these terms $G/B$ is associated to the empty set of nodes
and $G/P_{\min}$ is associated to a single node. It is also known that
such a presentation is unique, that is if $G_1/B_1\iso G_2/B_2$ where
$G$'s are simply connected semisimple and $B$'s are Borel then we can identify
$G_1=G_2$ and under this identification $B$'s are conjugate.

Now the proof of \ref{rigidity} goes as follows.
\begin{proof} Since $-K_{\cX/\Delta}$ is $\pi$-ample it follows that
  $X_0$ is Fano. By the main theorem of \cite{BLMS} and
  Kawamata-Shokurov contraction theorem every elementary contraction
  $\varphi_0: X_0\ra Y_0$ extends to an elementary contraction
  $\varphi: \cX\ra\cY$ relative over $\Delta$ with $\varphi_t: X_t\ra
  Y_t$ being an elementary contraction for every $t$. By our
  assumption, for $t\ne 0$ we have $X_t\iso G/B$, and by what we have
  said above the resulting elementary contraction $\varphi_t: X_t\ra
  Y_t$ is a $\PP^1$-bundle over $Y_t\iso G/P_{\min}$. Now by theorem
  \ref{generation} for $t\ne 0$ the variety $Y_t$ satisfies the
  assumptions of proposition \ref{equidimensional}. Therefore, by
  corollary \ref{proj-bundle}, $\varphi_0: X_0\ra Y_0$ is a
  $\PP^1$-bundle and by \ref{charOSWW} the theorem \ref{rigidity}
  follows.
\end{proof}

\subsection{A lemma}

\begin{lemma}\label{Poincare}
  Let $X_0$ and $X_1$ be a smooth projective varieties which are
  homeomorphic $X_0\simeq X_1$.  Assume that:
  \begin{enumerate}
  \item there exists a fiber bundle structure $\varphi_1:X_1\to Y_1$
    to a smooth variety $Y_1$ such that the fundamental group of $Y_1$
    acts trivially on the cohomology of the fiber,
  \item there exists a morphism $\varphi_0:X_0\to Y_0$ to a possibly
    singular variety $Y_0$,
  \item for some $k\geq 1$ the elements of degree $\leq k$, that is in
    $\HH^{\leq k}(Y_1)$, generate the ring $\HH^*(Y_1)$
  \item under identification given by the homeomorphism $X_0\simeq
    X_1$ we have
  $$\varphi^*_1(\HH^{\leq k}(Y_1))
  \subseteq \varphi^*_0(\HH^{\leq k}(Y_0))$$
  \end{enumerate}
  Then the dimension of every fiber of $\varphi_0$ does not exceed
  $\dim X_1-\dim Y_1$.
\end{lemma}

\begin{remark}The assumption about the trivial action of the fundamental group of  $Y_1$
  can by replaced by any other assumption which will guarantee that
  $\HH^*(X_1)$ is a free module over $\HH^*(Y_1)$.
\end{remark}
\begin{proof}
  We let $n=\dim X_0=\dim X_1$ and $f=\dim X_1-\dim Y_1$. Let
  $Z\subset X_0$ be a fiber of $\varphi_0$ and $z=\dim Z$. By
  $[Z]\in\HH^{2n-2z}(X_0)$ we denote the cohomology class of $Z$. To
  get a contradiction we assume that $z>f$.

  First we claim that for every $\alpha\in \HH^{>0}(Y_0)$ we have
  $[Z]\smile\varphi_0^*(\alpha)=0$. Indeed, let us look at the diagram
  of morphisms
  $$\begin{matrix}&&i\\
    &Z&\hookrightarrow&X_0\\
    \varphi_Z&\big\downarrow&\phantom{\Big\downarrow}&\big\downarrow&\varphi_0\\
    &pt&\hookrightarrow&Y_0\\ &&j\end{matrix}$$
  We have
  $$[Z]\smile\varphi_0^*(\alpha)=i_*i^*\varphi_0^*(\alpha)=
  i_*\varphi_Z^*j^*(\alpha)\,.$$ By our assumption $\deg(\alpha)>0$. Thus
  restricted to the point it vanishes, i.e.  $j^*(\alpha)=0$ and our
  claim follows.

  Now we claim that for any class $\beta\in \HH^{2z}(X_1)$ we have
  $[Z]\smile \beta=0$. Since $\varphi_1$ is a fibration of smooth
  projective varieties, the Serre spectral sequence degenerates by
  \cite[Thm.~II.1.2]{Blanchard} or \cite[Prop.~2.1]{Deligne}.
  Therefore the cohomology $\HH^*(X_1)$ is a free module over
  $\HH^*(Y_1)$.  Let $\{\gamma_i\}_{i \in I}$ be a homogeneous basis,
  with the degrees of $\gamma_i$'s at most $2f$.  Write
  \begin{equation}\label{lehi}\beta=\sum_I  \varphi_1^*(\delta_i)\smile \gamma_i\end{equation}
  with $\delta_i\in \HH^*(Y_1)$.  The elements $\varphi_1^*(\delta_i)$
  belong to the subring generated by
  $$\varphi_1^*(\HH^{\leq k}(Y_1))=\varphi_0^*(\HH^{\leq k}(Y_0))
  \subseteq\varphi_0^*(\HH^*(Y_0))$$
  Note that in the decomposition (\ref{lehi}) the nonzero $\delta_i$'s are of
  degree at least $2(z-f)>0$.  Thus $[Z]\smile
  \varphi_1^*(\delta_i)=0$ by what we proved above.  Then
  $$[Z]\smile\beta=\sum_I[Z]\smile
  \varphi_1^*(\delta_i) \smile\gamma_i=0$$
  Finally, by Poincar\'e duality, we conclude that $[Z]=0$ in
  $\HH^*(X_1)=\HH^*(X_0)$. Hence, since $X_0$ is projective, a
  contradiction.
\end{proof}

\subsection{Proof of proposition \ref{equidimensional}}
We will show that in the situation of \ref{equidimensional} the
varieties $X_0$ and $X_t$, $t\ne 0$, satisfy the assumptions of lemma
\ref{Poincare}, with $X_1=X_t$ and $k=2$.

Firstly we note that by the theorem of Ehresmann the family $\cX$ is
topologically locally trivial.  If needed we can shrink the base
$\Delta$ of $\pi:\cX\ra\Delta\ni 0$ and identify topology of fibers
$X_t$, for $t\ne 0$ with that of $X_0$, that is $X_0\simeq X_t$.  In
fact, for the purpose of topological argument we may assume $\Delta$
to be an open (in the classical topology) neighbourhood of $0$ which
is contractible. Then for every $t\in\Delta$ the inclusion
$\iota_{X_t}:X_t\hookrightarrow \cX$ is a homotopy equivalence.
Therefore we have isomorphisms
$$\HH^*(X_0)\stackrel{\iota^*_{X_0}}\lla
\HH^*(\cX)\stackrel{\iota^*_{X_t}}\lra \HH^*(X_t)\,.$$

\begin{lemma}\label{equality}
  In the set up of proposition \ref{equidimensional} let us assume
  that $-K_{\cX/\Delta}$ is $\varphi$-ample (we do not assume points
  (2) and (3) of this proposition). Then, under the identification
  $\HH^*(X_t)\simeq \HH^*(X_0)$, for $i=1,\ 2$ we have
  $$\varphi_t^* (\HH^i(Y_t)) \subseteq \varphi_0^*(\HH^i(Y_0))$$
\end{lemma}

\begin{remark}The opposite inclusion always holds as explained in
  section \ref{uwagi}.  \end{remark}

\begin{proof}
  We use results which follow from known vanishings related to Mori
  contractions. Since $-K_{X_t/Y_t}$ is $\varphi_t$-ample it follows
  that for every $t$ we have $\varphi_t^*(\HH^1(Y_t))=\HH^1(X_t)$, see
  e.g.~\cite[Prop.~2.3]{Asian}.

  Let $\cC(\cX/\cY)\subseteq \cC(\cX)$ be the cone of (classes of)
  curves contracted by $\varphi$.  By e.g.~\cite[Thm.~2.4]{Asian}, if
  $\varphi_R: \cX\ra \cY_R$ is a contraction of a Mori extremal ray
  $R\subseteq\cC(\cX/\cY)$ then $\varphi_R^*(\HH^2(\cY_R)) \subseteq
  \HH^2(\cX)$ is annihilated (orthogonal in terms of $\smile$-product)
  by the class of any curve contracted by $\varphi_R$. In fact
  $\varphi_R^*(\HH^2(\cY_R))=R^\perp$.  Since $\cC(\cX/\cY)$ is
  generated by Mori extremal rays whose contractions factor $\varphi$
  we get
  $$\varphi^*(\HH^2(\cY))=\bigcap_{R\subseteq\cC(\cX/\cY)}R^\perp$$
  Now by \cite[Prop.1.3]{Duke} the locus of every extremal ray
  contraction of $\cX$ dominates $\Delta$. Thus for $t\ne 0$ we have
  the left-hand-side inclusion
  $$\varphi_t^*(\HH^2(Y_t))\subseteq\bigcap_{R\subseteq\cC(\cX/\cY)}R^\perp
  \ =\ \varphi^*(\HH^2(\cY)) \ =\ \varphi_0^*(\HH^2(Y_0))$$ The
  right-hand-side equality follows because
  $\iota_{Y_0}:Y_0\hookrightarrow \cY$ may be assumed a homotopy
  equivalence (after possibly shrinking $\cY\rightarrow\Delta$). This
  concludes the proof of the lemma.
\end{proof}
From the preceding discussion and lemma \ref{equality} it follows
that assumptions of \ref{equidimensional} imply the set up of lemma
\ref{Poincare} hence proposition \ref{equidimensional} follows.

To prove corollary \ref{proj-bundle} we note that if $\varphi_t$ is a
$\PP^r$-bundle then we can choose a divisor class $H\in Pic(\cX)$ such
that its restriction $H_t$ to $X_t$ for $t\ne 0$ is a relative
$\cO(1)$ bundle for the projective bundle $\varphi_t: X_t\ra Y_t$. By
the arguments which we used in the proof of lemma \ref{equality} we
have $Pic(\cX)/ \varphi^*(Pic(\cY)) \iso \ZZ\cdot[H]$, hence
$(r+1)H+K_{\cX/\cY} \in \varphi^*Pic(\cY) $.  Moreover by proposition
\ref{equidimensional} $\varphi$ has all fibers of dimension $r$.
Thus, in view of Fujita's theorem \cite[2.12]{Fujita} the morphism
$\varphi: \cX\ra \cY$ is a $\PP^r$-bundle and the same concerns
$\varphi_0$.

\subsection{Generalities on cohomology of $G/P$}
Before giving the proof of \ref{generation} let us recall basic facts
about cohomology ring of homogeneous spaces. Let $G$ be a semisimple
complex algebraic group. We fix a maximal torus and a Borel subgroup
$\TT<B<G$. Let $P<G$ be a parabolic subgroup containing $B$.

The roots of the Lie algebra $\mathfrak{g}$ belong to
$\mathfrak{t}^*_{\ZZ}=\Hom(\TT,\CC^*)$.  Parabolic subgroups of $G$
which contain the fixed Borel group correspond to the subsets of the
set of simple roots. The empty set corresponds to the Borel group.  A
minimal parabolic subgroup $P=P_{\min}$ corresponds to a single root
$\alpha$.

Let $\mathfrak{t}^*_\QQ=\mathfrak{t}^*_\ZZ\otimes\QQ$. The Weyl group
$W_G=W=N(\TT)/\TT$ acts on $\mathfrak{t}^*_\QQ$. We fix an invariant
scalar product e.g.~the negative of the Killing form. The Weyl group
can be identified with the group generated by reflections in
$\alpha^\perp$ where $\alpha$ are simple roots. The Weyl group $W_P$
of a parabolic subgroup $P$ is the subgroup of $W$ generated by the
simple roots defining $P$.

Recall that for any parabolic subgroup $P<G$ the homogeneous variety
$G/P$ admits a decomposition into Schubert cells, therefore
$\HH^*(G/P,\ZZ)$ is a free abelian group generated by algebraic
cycles.  Further on we consider cohomology with rational
coefficients. Let $S=Sym(\mathfrak{t}_\QQ^*)$ be the symmetric algebra
generated by $\mathfrak{t}_\QQ^*$ and $S^W$ (resp. $S^{W_P}$) be the
subalgebra of $W$--invariants (resp $W_P$--invariants). By $S^W_+\subset S^W$ we denote the ideal of elements having positive degrees. The rational
cohomology ring of $G/P$ was computed in \cite[Th. 20.6(b)]{Borel} and
\cite[Thm.~5.5]{BGG}.
\begin{thm}[\cite{Borel, BGG}] With the notation as above we have
\begin{equation}\label{wzor}\HH^*(G/P)\simeq S^{W_P}/(S^W_+)\,,\end{equation}
where $(S^W_+)$ is the ideal in $S^{W_P}$ generated by $S^W_+$.
\end{thm}

The isomorphism in (\ref{wzor}) preserves the gradation under the
convention that the linear forms $\mathfrak{t}_\QQ^*\subset S$ live in
the second gradation. We sketch a short proof of this formula.

\begin{proof}
  Let $\EE_G$ be a contractible space with a free $G$--action.  We
  have a fibration
  \begin{equation}\label{rozwl}G/P\stackrel{\iota}\hookrightarrow
    \BB_P\stackrel{p}\twoheadrightarrow\BB_G\,.\end{equation}
  Here $\BB_P=\EE_G/P$, $\BB_G=\EE_G/G$ are the classifying spaces of
  Lie groups.  The cohomology algebra of $\BB_G$ was computed in
  \cite[Th. 20.3]{Borel}:
  \begin{equation}\label{HBG}\HH^*(\BB_G)\simeq S^W,\end{equation}
  Similarly
  \begin{equation}\label{HBP}\HH^*(\BB_P)\simeq S^{W_P}.\end{equation}
  The base and the fiber of (\ref{rozwl}) have cohomology concentrated
  in even degrees and therefore the Serre spectral sequence of the
  fibration degenerates. We obtain a surjection of algebras
  $$\HH^*(\BB_P)\twoheadrightarrow\HH^*(G/P)$$
  induced by $\iota$.  The image of $p^*(\HH^{>0}(\BB_G))$ vanishes
  when restricted to the fiber. The resulting surjection
  \begin{equation}\label{homologicznie}
    \HH^*(\BB_P)/(\HH^{>0}(\BB_G))\twoheadrightarrow \HH^*(G/P)
  \end{equation}
  of algebras is an isomorphism by counting the dimensions in each
  gradation.  Combining (\ref{homologicznie}) with (\ref{HBG}) and
  (\ref{HBP}) we obtain the formula (\ref{wzor}).
\end{proof}

In particular, if $P=B$ is the Borel subgroup, then
$\HH^{2*}(G/B)\simeq S/(S^W_+)$ is generated by $\HH^2(G/B)\simeq
\mathfrak{t}_\QQ^*$. An integral character $\chi\in
\Hom(\TT,\CC^*)=\mathfrak{t}_\ZZ^*$ corresponds to the first Chern
class of the line bundle defined by this character,
\cite[Cor.~4]{Gro}. The integral algebra structure $\HH^{*}(G/B,\ZZ)$
is much more complicated, see e.g. \cite{Toda}.

\subsection{Proof of theorem \ref{generation}}
If $P$ is a minimal parabolic subgroup corresponding to a simple root
$\alpha$ then $W_P=\ZZ/2$ is generated by the reflection in
$\alpha^\perp$. Therefore $$S^{W_P}=Sym(\alpha^\perp) \otimes
\QQ[\alpha^2]$$ Hence $S^{W_P}$ is generated by the linear and
quadratic forms.

Our argument is based on the formula \ref{wzor}. We will show that the
quotient $S^{W_P}/(S^W_+)$ is generated by the linear forms. It is
enough to show that $\alpha^2$ can be expressed by $W_P$-invariant
linear forms modulo a quadratic $W$-invariant form. Let $Q$ be the
quadratic form corresponding to the invariant scalar product.

We choose vectors $\beta_1,\beta_2,\dots,\beta_{k}$ which form an
orthogonal basis of $\alpha^\perp$. Then for some nonzero numbers
$a,b_1,b_2,\dots b_k$ the quadratic form $Q$ can be written as
$$Q=a \alpha^2+\sum_{i=1}^kb_i\beta_i^2\,.$$
Thus
$$ \alpha^2\equiv -\sum_{i=1}^k\frac{b_i}{a}\beta_i^2\;\mod\;S^W\,.$$

This concludes the proof of theorem \ref{generation}.

\section{Remarks}\label{uwagi}

\begin{rem}
  We note that the assumption that $-K_{\cX/\Delta}$ is $\pi$-ample in
  Theorem \ref{rigidity} is necessary because
  $X_t\iso\PP^1\times\PP^1=\PP_{\PP^1}(\cO\oplus\cO)$ can be
  specialized to $X_0\iso\PP_{\PP^1}(\cO(-a)\oplus\cO(a))$ being an
  arbitrary even Hirzebruch surface. However it make sense to ask if
  the above rigidity theorem remains true remains true for families of
  projective manifolds (that is: $-K_{\cX/\Delta}$ not necessarily assumed
  to be $\pi$-ample) if one assumes that $G$ is simple (that is: with
  $X_t$ irreducible complete flag variety).
\end{rem}

Also, a naive extension of Theorem \ref{rigidity} to non-complete
flags $G/P$, where is $P<B$ is any parabolic, does not work.

\begin{ex} The partial flag variety $X=Fl_{1,2}(V)$ of lines
  and planes in $V=\CC^{2n+2}$ admits a fibration to the Grassmannian of
  planes $Grass_2(V)$ and to $\PP^{2n+1}$. Let
  $\Omega=\Omega^1_{\PP^{2n+1}}$ be the cotangent bundle to
  $\PP^{2n+1}$. One can identify $X$ with the Grothendieck
  projectivization of $\Omega(2)=\Omega\otimes\cO(2)$. Let
  $\omega\in\bigwedge^2V$ be a symplectic form. We can use the
  second exterior power of the dual Euler sequence
  \begin{equation}\label{Euler}0\longrightarrow \Omega^2(2)
    \longrightarrow \bigwedge^2V
    \longrightarrow\Omega(2)\longrightarrow 0\end{equation}
  to identify $\omega$ with a section of $\Omega(2)$ and thus obtain a surjective map
  $\cT_{\PP^{2n+1}}(-1)\ra\cO(1)$. The kernel $\cN$ is called the
  null-correlation bundle; see \cite[Sect.~4.2]{OSS} where one can
  find details of this construction. In other terms $\cN(1)\subset\cT_{\PP^{2n+1}}$ is the contact distribution associated with the symplectic form $\omega$.

  The (Grothendieck) projectivization $\PP(\cN^*(1)) \subset
  \PP(\Omega(2))$ is the incidence variety of lines and
  $\omega$-isotropic planes.  The extension
  \begin{equation}\label{extension}0\lra\cO\lra
    \Omega(2)\lra\cN^*(1)\lra 0\end{equation}
  can be specialized to $\cN^*(1)\oplus\cO$ and the projectivization
  $X_t=\PP(\Omega(2))$ specializes to $X_0=\PP(\cN^*(1)\oplus\cO)$.
  Consequently, the $\PP^1$-bundle
  \begin{equation}\label{wiazka}Fl_{1,2}(V)=\PP(\Omega(2))\ra
    Grass_2(V)\end{equation} specializes to a map of
  $\PP(\cN\oplus\cO)$ to a cone over the isotropic Grassmannian. The
  latter morphism has the unique fiber $Z$ over the vertex of the cone
  which is not $\PP^1$ but $\PP^{2n+1}$.

  One easily checks that here the image of $\HH^*(Y_0)$ in
  $\HH^*(X_0)$ is smaller than the image of $\HH^*(Y_t)$
  e.g.~computing the ranks of these groups.  At this point the proof
  of proposition \ref{equidimensional} breaks in our example.

  Let us see what happens with the class $[Z]$.  The $\PP^1$-bundle
  (\ref{wiazka}) is defined by global sections of $\Omega(2)$. Thus
  the pull-back of the ample generator of Picard group of $Grass_2(V)$
  is the relative $\cO_{\PP(\Omega(2)/\PP^{2n+1}}(1)$ on
  $\PP(\Omega(2))$.  We denote it by $D$.  On the other hand by $H$ we
  denote the hyperplane class from $\PP^{2n+1}$; it is the relative
  $\cO(1)$ on the projectivization of the universal bundle over
  $Grass_2(V)$ which makes the incidence variety $Fl_{1,2}(V)$.  Let
  $\sum_{i=0}^{2n+1}c_iH^i$, with $c_i\in\ZZ$ be the Chern class of
  $\Omega(2)$. The numbers $c_i$ can be calculated from the Euler
  sequence.  By Leray-Hirsch formula we have
  \begin{equation}\label{leray-hirsch}
  \HH^{2*}(\PP(\Omega(2))\simeq\QQ[H,D]/(H^{2n+2},f(H,D))
  \end{equation}
  where the polynomial $f$ is given by the
  formula:
  \begin{equation}\label{chern}
  f(H,D)=\sum_{i=0}^{2n+1}(-1)^ic_iH^iD^{2n+1-i}\,.  \end{equation}
  Note that $c_{2n+1}=0$ since the bundle $\Omega(2)$ has a nowhere
  vanishing section. We can write \begin{equation}\label{splitting2}
  f(H,D)=f_0(H,D)\cdot D \end{equation} for some homogeneous
  polynomial $f_0\in\ZZ[H,D]$ of degree $2n$.  Comparing the formula
  (\ref{splitting2}) with the splitting sequence (\ref{extension}) we
  see that class of the fiber $Z\iso \PP^{2n+1}$ in
  $\HH^{4n}(\PP(\cN\oplus\cO))$ via the isomorphism
  (\ref{leray-hirsch}) it can be identified with $f_0(H,D)$. In
  particular it is annihilated by the class $D$. Since $H$ is the
  class of the relative $\cO(1)$ on $Fl_{1,2}(V)$ treated as $\PP^1$
  over $Grass_2(V)$ we can write
  $$[Z]=\varphi_t^*(\alpha)+\varphi_t^*(\beta)\smile H\,,$$
  for some classes $\alpha\in \HH^{4n}(Grass_2(\CC^{2n+2}))$ and
  $\beta\in \HH^{4n-2}(Grass_2(\CC^{2n+2}))$. Since $[Z]\smile D=0$,
  the classes $\alpha$ and $\beta$ are annihilated by
  $D$. Applying Hard Lefschetz theorem we conclude that
  $\alpha$ belongs to the primitive cohomology and $\beta=0$.  In
  terms of Schubert classes $\alpha$ is proportional to
  $\sum_{i=0}^{n}(-1)^iS_{[n+i,n-i]}$. Here $S_{[n+i,n-i]}$ is the cohomology class of the Schubert variety defined by the partition $[n+i,n-i]$. This shows that
  $[Z]=\varphi_t^*(\alpha)$ and $\alpha$ is not effective.
\end{ex}

In the example above we observe that the cohomology of $Y_0$ can be
smaller than that of a general fiber. In fact
$$\varphi_0^* (\HH^*(Y_0))\subseteq \varphi_t^*(\HH^*(Y_t))\,.$$
  This is a general rule, but it rather concerns the image of
  $\HH^*(Y_0)$ in $\HH^*(X_0)$. Suppose $\Delta$ is a disc in $\CC$. We have a commutative diagram:
 \begin{equation}\label{diagram}\begin{matrix}
 \HH^*(X_0)&\stackrel{\simeq}\lla&H^*(\cX)&\stackrel{\simeq}\lra&H^*(X_t)\\
 \varphi_0^*\Big\uparrow\phantom{\varphi_0^*}&\phantom{\Bigg\uparrow}&\varphi^*\Big\uparrow\phantom{\varphi^*}&&\varphi_t^*\Big\uparrow&\phantom{\varphi_t^*}\\
 \HH^*(Y_0)&\stackrel{\simeq}\lla&\HH^*(\cY)&\lra&\HH^*(Y_t)\\
\end{matrix}\end{equation}
 In general, for fibrations over
$\Delta\setminus\{0\}$ there is a monodromy operator (usually denoted
by $T$) given by the action of the generator of
$\pi_1(\Delta\setminus\{0\},t)\simeq\ZZ$ on $\HH^*(Y_t)$. The image
of $$\HH^*(Y_0)\ra\HH^*(Y_t)$$ is contained in the invariant subspace
of monodromy.  By the ``invariant cycle theorem'' \cite{clemens},
\cite{morrison} if $\cY$ is a K\"ahler manifold then we have
an equality
$${\rm im}\left(\HH^*(Y_0)\ra\HH^*(Y_t)\right)=\HH^*(Y_t)^T\,.$$
For the maps $\cY\to\Delta$ obtained by contraction of a smooth map
$\varphi: \cX\to\Delta$ the monodromy is trivial, since the cohomology
of $Y_t$ embeds into the cohomology of $X_t$ where clearly the
monodromy is trivial. Hence
$\varphi_0(H^*(Y_0))=\varphi_t(H^*(Y_t))$. Therefore by lemma
\ref{Poincare} we get the following.

\begin{cor}With the notation of proposition \ref{equidimensional}
  assume (2) and that $\cY$ is smooth.  Then
  the dimension of every fiber of $\varphi_0$ does not exceed $\dim
  X_0-\dim Y_0$.
\end{cor}

We conclude with the remark that further understanding deformation in the context of Mori contraction
undoubtedly leads to study the vanishing cycle sheaf on the contracted
space $\cY$.


\begin{thebibliography}{10}

\bibitem{BGG}
I.~N. Bern{\v{s}}te{\u\i}n, I.~M. Gel{\cprime}fand, and S.~I. Gel{\cprime}fand.
\newblock Schubert cells, and the cohomology of the spaces {$G/P$}.
\newblock {\em Uspehi Mat. Nauk}, 28(3(171)):3--26, 1973.

\bibitem{Blanchard}
Andr{\'e} Blanchard.
\newblock Sur les vari\'et\'es analytiques complexes.
\newblock {\em Ann. Sci. Ecole Norm. Sup. (3)}, 73:157--202, 1956.

\bibitem{Borel}
Armand Borel.
\newblock {\em Topics in the homology theory of fibre bundles}, volume 1954 of
  {\em Lectures given at the University of Chicago}.
\newblock Springer-Verlag, Berlin-New York, 1967.

\bibitem{clemens}
C.~H. Clemens.
\newblock Degeneration of {K}\"ahler manifolds.
\newblock {\em Duke Math. J.}, 44(2):215--290, 1977.

\bibitem{Deligne}
P.~Deligne.
\newblock Th\'eor\`eme de {L}efschetz et crit\`eres de d\'eg\'en\'erescence de
  suites spectrales.
\newblock {\em Inst. Hautes \'Etudes Sci. Publ. Math.}, (35):259--278, 1968.

\bibitem{Fujita}
Takao Fujita.
\newblock On polarized manifolds whose adjoint bundles are not semipositive.
\newblock In {\em Algebraic geometry, {S}endai, 1985}, volume~10 of {\em Adv.
  Stud. Pure Math.}, pages 167--178. North-Holland, Amsterdam, 1987.

\bibitem{Gro}
Alexander Grothendieck.
\newblock Torsion homologique et sections rationelles.
\newblock In {\em S\'eminaire {C}. {C}hevalley; 2e ann\'ee: 1958. {A}nneaux de
  {C}how et applications}, volume~3, pages exp no 5, 1--29. Secr\'etariat
  math\'ematique, 11 rue Pierre Curie, Paris, 1958.

\bibitem{HwangMokInvMath}
Jun-Muk Hwang and Ngaiming Mok.
\newblock Rigidity of irreducible {H}ermitian symmetric spaces of the compact
  type under {K}\"ahler deformation.
\newblock {\em Invent. Math.}, 131(2):393--418, 1998.

\bibitem{HwangMokAnnSciENS}
Jun-Muk Hwang and Ngaiming Mok.
\newblock Deformation rigidity of the rational homogeneous space associated to
  a long simple root.
\newblock {\em Ann. Sci. \'Ecole Norm. Sup. (4)}, 35(2):173--184, 2002.

\bibitem{KollarMori}
J{\'a}nos Koll{\'a}r and Shigefumi Mori.
\newblock {\em Birational geometry of algebraic varieties}, volume 134 of {\em
  Cambridge Tracts in Mathematics}.
\newblock Cambridge University Press, Cambridge, 1998.
\newblock With the collaboration of C. H. Clemens and A. Corti, Translated from
  the 1998 Japanese original.

\bibitem{morrison}
David~R. Morrison.
\newblock The {C}lemens-{S}chmid exact sequence and applications.
\newblock In {\em Topics in transcendental algebraic geometry ({P}rinceton,
  {N}.{J}., 1981/1982)}, volume 106 of {\em Ann. of Math. Stud.}, pages
  101--119. Princeton Univ. Press, Princeton, NJ, 1984.

\bibitem{OSWW}
Gianluca Occhetta, Luis~E. Sol{\'a}~Conde, Kiwamu Watanabe, and Jaros{\l}aw~A.
  Wi{\'s}niewski.
\newblock {F}ano manifolds whose elementary contractions are smooth $\mathbb
  {P}^1$-fibrations.
\newblock {to appear in  Annali della Scuola Normale Superiore di Pisa,
  Classe di Scienze}, 
\newblock Preprint arXiv:{\tt 1407.3658, DOI:~10.2422/2036-2145.201508\_007}.

\bibitem{OSS}
Christian Okonek, Michael Schneider, and Heinz Spindler.
\newblock {\em Vector bundles on complex projective spaces}.
\newblock Modern Birkh\"auser Classics. Birkh\"auser/Springer Basel AG, Basel,
  2011.
\newblock Corrected reprint of the 1988 edition, With an appendix by S. I.
  Gelfand.

\bibitem{PasquierPerrin}
Boris Pasquier and Nicolas Perrin.
\newblock Local rigidity of quasi-regular varieties.
\newblock {\em Math. Z.}, 265(3):589--600, 2010.

\bibitem{SiuCrelle}
Yum~Tong Siu.
\newblock Errata: ``{N}ondeformability of the complex projective space'' [{J}.
  {R}eine {A}ngew.\ {M}ath.\ {\bf 399} (1989), 208--219; {MR}1004139
  (90h:32048)].
\newblock {\em J. Reine Angew. Math.}, 431:65--74, 1992.

\bibitem{Toda}
Hirosi Toda.
\newblock On the cohomology ring of some homogeneous spaces.
\newblock {\em J. Math. Kyoto Univ.}, 15:185--199, 1975.

\bibitem{Duke}
Jaros{\l}aw~A. Wi{\'s}niewski.
\newblock On deformation of nef values.
\newblock {\em Duke Math. J.}, 64(2):325--332, 1991.

\bibitem{Asian}
Jaros{\l}aw~A. Wi{\'s}niewski.
\newblock Cohomological invariants of complex manifolds coming from extremal
  rays.
\newblock {\em Asian J. Math.}, 2(2):289--301, 1998.

\bibitem{BLMS}
Jaros{\l}aw~A. Wi{\'s}niewski.
\newblock Rigidity of the {M}ori cone for {F}ano manifolds.
\newblock {\em Bull. Lond. Math. Soc.}, 41(5):779--781, 2009.

\end{thebibliography}
\bibliographystyle{plain}
\def\cprime{$'$}

\end{document}